\documentstyle[amscd,amssymb,verbatim,12pt]{amsart}
\newcommand{\area}{\mbox{\rm area}}

\newcommand{\bb}{\mbox{\rm b}}

\newcommand{\BE}{\widetilde{\mbox{\rm Ric}}_\infty}
\newcommand{\BEq}{\widetilde{\mbox{\rm Ric}}_q}

\newcommand{\cd}{\mbox{\rm CD}}

\newcommand{\const}{\mbox{\rm const.}}

\newcommand{\diam}{\mbox{\rm diam}}

\newcommand{\grad}{\mbox{\rm grad}}

\newcommand{\Hess}{\mbox{\rm Hess}}

\newcommand{\HH}{\mbox{\rm H}}

\newcommand{\R}{{\Bbb R}}

\newcommand{\Ric}{\mbox{\rm Ric}}
\newcommand{\Riem}{\mbox{\rm Riem}}

\newcommand{\Tr}{\mbox{\rm Tr}}

\newcommand{\vol}{\mbox{\rm vol}}
\newcommand{\Z}{{\Bbb Z}}
\theoremstyle{plain}
\newtheorem{definition}{Definition}

\newtheorem{lemma}{Lemma}
\newtheorem{theorem}{Theorem}

\numberwithin{equation}{section}
\renewcommand{\rm}{\normalshape}
\begin{document}
\title{Some Geometric Properties of the Bakry-\'Emery-Ricci Tensor}
\author{John Lott}
\address{Department of Mathematics\\
University of Michigan\\
Ann Arbor, MI  48109-1109\\
USA}
\email{lott@@umich.edu}
\thanks{Research supported by NSF grant DMS-0072154}
\date{October 16, 2002}
\maketitle
\begin{abstract}
The Bakry-\'Emery tensor gives an analog of the Ricci tensor for a
Riemannian manifold with a smooth measure.  We show that some of
the topological consequences of having a positive or nonnegative
Ricci tensor are also valid for the Bakry-\'Emery tensor. We show
that the Bakry-\'Emery tensor is nondecreasing under a Riemannian
submersion whose fiber transport preserves measures up to constants.
We give some
relations between the Bakry-\'Emery tensor and measured 
Gromov-Hausdorff limits.
\end{abstract}

\section{Introduction} 

When considering the metric
structure of manifolds with lower Ricci curvature bounds, it is natural
to carry along the extra structure of a measure and consider metric-measure
spaces.  This is especially relevant for collapsing, and
has been discussed by Cheeger-Colding
\cite{Cheeger-Colding (1997),Cheeger-Colding (2000a),Cheeger-Colding (2000b)},
Fukaya \cite{Fukaya (1987)}
and Gromov \cite[Chapter $3 \frac12$]{Gromov (1999)}.

In this paper we consider smooth metric-measure spaces.
Let $M$ be a connected
$n$-dimensional Riemannian manifold, with metric $g$. Let
$d\vol_M$ denote the Riemannian density on $M$.  Let $\phi$ be a smooth
positive function on $M$.  Then $(M, \phi \: d\vol_M)$ is
a smooth metric-measure space. For reasons coming from the study of
diffusion processes, Bakry and \'Emery
\cite{Bakry-Emery (1985)} defined a generalization of
the Ricci tensor of $M$ by
\begin{equation} \label{1.1}
\BE \: = \: \Ric \: - \: \Hess (\ln \phi).
\end{equation}
In terms of indices,
$\left( \BE \right)_{\alpha \beta} \: = \: \Ric_{\alpha \beta} \: - \: 
(\ln \phi)_{; \alpha \beta}$.

It turns out that the Bakry-\'Emery tensor (\ref{1.1})
has interesting connections
to logarithmic Sobolev inequalities, isoperimetric inequalities
and heat semigroups.  We refer to
\cite{Ane (2000)} and
\cite{Ledoux (2000)} for information on these connections.
(In fact, Bakry and \'Emery defined their tensor in a more abstract
setting than what we consider.)

We are interested in the geometric implications of bounds on the
Bakry-\'Emery tensor. As in \cite{Qian (1997)}, let us define
a related tensor $\BEq$. Given $q \in (0, \infty)$, put
\begin{align} \label{1.2}
\widetilde{\Ric}_q \: & = \:
\Ric \: - \: 
\Hess(\ln \phi) \: - \:
\frac{1}{q} \: d \ln \phi \:  \otimes \: d \ln \phi \\
& = \: \Ric \: - \: 
\frac{\Hess(\phi)}{\phi} \: + \:
\left( 1 \: - \: \frac{1}{q} \right) \: \frac{d \phi}{\phi} \:  
\otimes \: \frac{d \phi}{\phi} \notag \\
& = \:
\Ric \: - \: q \:
\frac{\Hess \left( \phi^{\frac{1}{q}} \right)}{\phi^{\frac{1}{q}}}. \notag
\end{align}
Clearly, if $\BEq \: \ge \: r \: g$ then
$\BE \: \ge \: r \: g$. 
In the terminology of \cite{Bakry (1994)}, a condition of the form
$\BEq \: \ge \: r \: g$ implies a curvature-dimension inequality
$\cd(r, n+q)$.

Our first result extends some classical topological  results about the Ricci 
tensor (i.e. when $\phi$ is constant)
to the setting of the Bakry-\'Emery tensor.

\begin{theorem} \label{theorem1} Suppose that $M$ is connected and closed. \\
1. If $\BE \: > 0$ then $\pi_1(M)$ is finite. \\
2. If $\BEq \: \ge \: 0$ and $q \in (0, \infty)$ 
then $\pi_1(M)$ has a finite-index free abelian
subgroup of rank at most $n$.\\
3.  If $\BE \: \ge 0$ then $\HH^1(M; \R)$ is isomorphic to the
linear space of parallel $1$-forms on $M$ whose pairing with
$\grad (\phi)$ vanishes identically.  In particular, if
$\BE \: \ge 0$ then  $\bb_1(M)  \: \le \: n$.
If $\BE \: \ge 0$ and $\bb_1(M) \: = \: n$ then 
$M$ is a flat torus and $\phi$ is
constant. \\
4. If $\BE \: < \: 0$ then the isometry group of $(M, g)$ is
finite. \\
5. If $\BE \: \le \: 0$ then any Killing vector field on $(M, g)$
is parallel and annihilates $\phi$.
\end{theorem}
\noindent
{\bf Remark : } Theorem \ref{theorem1}.2 is a strengthening of
\cite[Theorem 6]{Qian (1997)}, which says that if 
$\BEq \: \ge \: 0$ and $q \in (0, \infty)$ 
then $\pi_1(M)$ has polynomial growth of order at most $n \: + \: q$. \\

The proofs of parts 3-5 of Theorem \ref{theorem1} use a Bochner-type
identity.  If the pair $(g, \phi)$ is only $C^0 \cap H^1$-regular then
one can use this identity to still make sense of the notion
$\BE \: \ge \: r \: g$ or $\BEq \: \ge \: r \: g$
(see Definition \ref{definition1} of Section \ref{Section 2}).

An important result in the study of manifolds of nonnegative
sectional curvature is O'Neill's theorem, which says that
sectional curvature is nondecreasing under a Riemannian
submersion
\cite[Chapter 9]{Besse (1987)}. We show that there is a Ricci
analog of O'Neill's theorem, provided that one uses the Bakry-\'Emery
tensor and assumes that the fiber transport of the Riemannian submersion
preserves measures up to multiplicative constants.

Suppose that a Riemannian submersion  $p \: : \: M \rightarrow B$
has compact fiber $F$. Put $F_b \: = \: p^{-1}(b)$. Given
a smooth curve $\gamma \: : \: [0,1] \rightarrow B$ and a point
$m \in F_{\gamma(0)}$, let $\rho(m)$ be the endpoint 
$\overline{\gamma}(1)$ of the horizontal lift
$\overline{\gamma}$ of $\gamma$ that
starts at $\overline{\gamma}(0)$. Then $\rho$ is the
fiber transport diffeomorphism
from $F_{\gamma(0)}$ to $F_{\gamma(1)}$. 

Given the positive function $\phi^M$ on $M$, define
$\phi^B$, a smooth positive function on $B$, by
\begin{equation} \label{1.3} 
p_*(\phi^M \: d\vol_M) \: = \: \phi^B \: d\vol_B.
\end{equation}
Let $\BE^M$ and $\BE^B$ denote the corresponding Bakry-\'Emery tensors.
Let $d\vol_F$ denote the fiberwise Riemannian density.

\begin{theorem} \label{theorem2}
Suppose that fiber transport 
preserves the fiberwise measure $\phi_M \: d\vol_F$ up to a 
multiplicative constant, 
i.e. for any smooth curve 
$\gamma \: : \: [0,1] \rightarrow B$,
there is a constant $c_\gamma \: > \: 0$ such that
$\rho^* \left( \phi^M \Big|_{F_{\gamma(1)}} \: d\vol_{F_{\gamma(1)}} \right) 
\: = \: c_\gamma \:
\phi^M \Big|_{F_{\gamma(0)}} \: d\vol_{F_{\gamma(0)}}$. \\
1. For any $r \in \R$, if $\BE^M \: \ge \: r \: g^M$ then
$\BE^B \: \ge \: r \: g^B$. \\
2. Suppose in addition that $\phi^M \: = \: 1$. Put $q \: = \: \dim(F)$.
For any $r \in \R$, if $\Ric^M \: \ge \: r \: g^M$ then
$\BEq^B \: \ge \: r \: g^B$. \\
\end{theorem}

Using Theorem \ref{theorem2}, we show a relationship between
$\BEq$ and collapsing.

\begin{theorem} \label{theorem3}
1. Given $r \in \R$ and an integer $q \: \ge \: 2$, let $(B, \phi)$ be
a smooth closed
measured Riemannian manifold with $\BEq^B \: \ge \: r \: g^B$. Then
$(B, \phi)$ is the measured Gromov-Hausdorff limit of a sequence of
$(n \: + \: q)$-dimensional closed Riemannian manifolds $(M_i, g_i)$ with
$\Ric(M_i, g_i) \: \ge \: r \: g_i$.\\
2. Let $\{(M_i, g_i)\}_{i=1}^\infty$ be a sequence of 
$N$-dimensional
connected closed Riemannian manifolds
with sectional curvatures bounded above in absolute value by $\Lambda$
and diameters bounded above by $D$, for some $D, \Lambda \in \R^+$.
Let $(X, \mu)$ be a limit point for 
$\{(M_i, g_i)\}_{i=1}^\infty$ in the measured Gromov-Hausdorff topology.  
Suppose that for some $r \in \R$ and all $i \in \Z^+$, 
$\Ric(M_i, g_i) \: \ge \: r \: g_i$. 
Suppose that $X$ is an $n$-dimensional closed manifold. 
Put $q \: = \: N \: - \: n$.\\
a. If $q \: = \: 0$
then $X$ has $\Ric \: \ge \: r g$ in the generalized sense of Definition
\ref{definition1} below.\\
b. If $q \: > \: 0$ then
$X$ has $\BEq \: \ge \: r \: g$ in the generalized sense
of Definition \ref{definition1} below.
\end{theorem}

Finally, we give a condition in terms of distances and masses that is
equivalent to having Bakry-\'Emery tensor bounded below by $r$.
If ${\cal O}$ is a measurable subset of $M$, put
\begin{equation} \label{1.4}
\vol_\phi({\cal O}) \: = \:
\int_{{\cal O}} \phi \: d\vol_M.
\end{equation}
Following \cite[Section 5.45]{Gromov (1999)}, we define the notion of
a distance tube in $M$. Let $T_0$ be a closed subset of $M$.  A subset
$T \subset M$ containing $T_0$ 
is a distance tube with base $T_0$ if for all $t \in T$,
there is a segment $s \subset T$ from some $t_0 \in T_0$ to $t$ with
length $l(s) \: = \: d(t, T_0)$.
For $0 \: < \: u_1 \: < \: u_2$, define the distance annulus
\begin{equation} \label{1.5}
A(u_1, u_2) \: = \: \{t \in T \: : \: u_1 \: \le \: d(t, T_0) \: \le \:
u_2\}.
\end{equation}
Given $c \in \R$, put
\begin{equation}
\widehat{v}(u_1, u_2, c) \: = \: 
\int_{u_1}^{u_2}
e^{- \: \frac{r}{2} x^2 \: + \: c \: x} \: dx.
\end{equation}

\begin{theorem} \label{theorem4}
Suppose that $\BE(M, g, \phi) \: \ge \: r \: g$ for some $r \in \R$.
Given numbers $0 \: < \: u_1 \: < \: u_2 \: < \: u_3$, we assume that
the tube $T$ is a disjoint union of segments $s$, starting at $T_0$,
of length at least $u_3$. We also assume that
$\vol_\phi(A(u_2, u_3)) \: > \: 0$.
Suppose that for some $c \in \R$,
\begin{equation} \label{1.6}
\frac{\vol_\phi(A(u_2, u_3))}{\vol_\phi(A(u_1, u_2))} \: \le \:
\frac{\widehat{v}(u_2, u_3, c)}{\widehat{v}(u_1, u_2, c)}.
\end{equation}
Then there is a subtube
$T^\prime \subset T$ consisting of a union of segments $s$ from $T_0$, 
such that \\
1. \begin{equation} \label{1.7}
\frac{\vol_\phi(T^\prime \cap A(u_1, u_2))}{\vol_\phi(A(u_1, u_2))}
 \: \ge \: 1 \: - \:
\frac{\vol_\phi(A(u_2, u_3))}{\vol_\phi(A(u_1, u_2))} \:
\left( 
\frac{\widehat{v}(u_2, u_3, c)}{\widehat{v}(u_1, u_2, c)} \right)^{-1},
\end{equation}
2. If a segment $s \subset T$, starting from $T_0$,
intersects $T^\prime \cap A(u_2, u_3)$
then $s \subset T^\prime$, and\\
3. For all $u_4 \: > \: u_3$,
\begin{equation} \label{1.8}
\frac{\vol_\phi(T^\prime \cap A(u_3, u_4))}{\vol_\phi(
T^\prime \cap A(u_2, u_3))} \: \le \:
\frac{\widehat{v}(u_3, u_4, c)}{\widehat{v}(u_2, u_3, c)}.
\end{equation}

Conversely, suppose that there is a number $r \in \R$ so that
for each tube $T$ and $c \in \R$ satisfying (\ref{1.6}), there is a
subtube $T^\prime$ with the above properties. Then
$\BE(M, g, \phi) \: \ge \: r \: g$.
\end{theorem}

In Sections \ref{Section 2}-\ref{Section 5}  
we prove Theorems \ref{theorem1}-\ref{theorem4}, respectively.
In Section \ref{Section 6} we make some remarks.

I thank Max Karoubi for his hospitality at the Universit\'e de Paris VII,
and Thierry Coulhon and Sasha Grigor'yan for their hospitality at the
Institut Henri Poincar\'e, while part of this research was performed.

\section{Proof of Theorem \ref{theorem1}}  \label{Section 2}
We first prove parts 1 and 2 of the theorem.
If $\BE > 0$ then $\BEq > 0$ for some $q \in (0, \infty)$.
Increasing $q$ if necessary,  we may assume without 
loss of generality that $q$ is an integer greater than one. 
Thus for parts 1 and 2, it is enough to consider the case when
$\BEq > 0$ or $\BEq \ge 0$, for some integer
$q$ greater than one.
 
Given $i \in \Z^+$,
consider $S^q \times {M}$ with the warped product metric
$g^{S^q \times {M}} \: = \: 
g^{{M}} \: + \: i^{-2} \:
{\phi}^{\frac{2}{q}} \: 
g^{S^q}$.
Let $p \: : \: S^q \times {M} \rightarrow 
{M}$ be the projection. 
Let $\overline{X}$ be the horizontal lift 
to $S^q \times {M}$ of a vector field $X$ on 
${M}$ and let $\overline{U}$ be a vertical vector field on 
$S^q \times {M}$.
From 
\cite[Proposition 9.106]{Besse (1987)},
\begin{align} \label{2.14}
\Ric^{S^q \times {M}}(\overline{X}, \overline{X}) & = 
p^* \left( \Ric^{{M}}({X}, {X}) \: - \: q \: 
\frac{\Hess({\phi}^{\frac{1}{q}})({X}, 
{X})}{{\phi}^{\frac{1}{q}}} \right), \\
\Ric^{S^q \times {M}}(\overline{X}, \overline{U}) & = 
0 \notag \\
\Ric^{S^q \times {M}}(\overline{U}, \overline{U}) & =
\Ric^{S^q}(\overline{U}, \overline{U}) \: + \:
(\overline{U}, \overline{U}) \:
p^*  \left( - \: \frac{\nabla^2 
{\phi}^{\frac{1}{q}}}{{\phi}^{\frac{1}{q}}} 
\: - \: (q \: - \: 1) \:
\frac{|\nabla 
{\phi}^{\frac{1}{q}}|^2}{{\phi}^{\frac{2}{q}}} \right). 
\notag
\end{align}
Taking
$i \rightarrow \infty$, we see that if $\BEq(M, g, \phi) \: \ge \:
r \: g$ then
$\left({M}, g^{{M}}, {\phi} \right)$ 
is the limit
of a sequence of $(n \: + \: q)$-dimensional manifolds with
Ricci curvature bounded below by $r$. If $r$ is positive then
from Myers' theorem, $\pi_1(S^q \times M) \cong \pi_1(M)$ is finite. 
This proves part 1 of the theorem.

Now suppose that $r \: \ge \: 0$. For $i$ large, the warped product metric
on $S^q \times M$ has nonnegative Ricci curvature.  There is a 
$k \ge 0$ so that 
$\pi_1(S^q \times M) \cong \pi_1(M)$ has a finite-index free abelian
subgroup of rank $k$ and the universal cover
$S^q \times \widetilde{M}$ has an isometric splitting as
$\R^k \times Y^{n+q-k}$, where $Y$ is closed and simply-connected
\cite{Cheeger-Gromoll (1971)}.
Considering the cohomology groups of
$S^q \times \widetilde{M} \: \cong \: \R^k \times Y^{n+q-k}$,
it follows that
\begin{equation}
q \: + \: \max \{j \: : \: \HH^j(\widetilde{M}; \Z) \neq 0 \} \: = \:
n+q-k.
\end{equation}
Then $k \: = \: n \: - \: 
\max \{j \: : \: \HH^j(\widetilde{M}; \Z) \neq 0 \} \: \le \: n$, which
proves part 2 of the theorem.

To prove the rest of the theorem,
if $V$ is a vector field on $M$, let $V^\sharp$ denote the dual $1$-form.
If $\omega$ is a $1$-form on $M$, let $\omega_\sharp$ denote the dual 
vector field.
Let $i_V$ denote interior multiplication with respect to $V$ and
let ${\cal L}_V$ denote Lie differentiation with respect to $V$.

If $T$ is a tensor field on $M$, let $(T, T) \in C^\infty(M)$ be
the inner product coming from the Riemannian metric $g$.
Put
\begin{equation} \label{2.1} 
\langle T, T \rangle \: = \: 
\int_M (T, T)(m) \: \phi(m) \: d\vol_M(m). 
\end{equation}
Let $(\Omega^*(M), d)$ denote the de Rham complex of $M$.
Let $\delta$ be the formal adjoint of $d$ with respect to the
Riemannian metric $g$, i.e. in the case $\phi \: = \: 1$, and let
$\widetilde{\delta}$ be the formal adjoint of $d$ with respect to
$\langle \cdot , \cdot \rangle$. Then
\begin{equation} \label{2.2} 
\widetilde{\delta} \: = \: \delta \: - \:
i_{(d \ln \phi)_\sharp}.
\end{equation}

Put $\triangle \: = \: 
d \delta \: + \: \delta d$ and 
$\widetilde{\triangle} \: = \: 
d \widetilde{\delta} \: + \: \widetilde{\delta} d$. Then
\begin{equation} \label{2.3} 
\widetilde{\triangle} \: = \: \triangle \: - \: 
d i_{(d \ln \phi)_\sharp} \: - \: i_{(d \ln \phi)_\sharp} d \: = \:
\triangle \: - \: {\cal L}_{(d \ln \phi)_\sharp}.
\end{equation}

The Bochner identity says that if $\omega$ is a
$1$-form then there is an equality of functions on $M$ :
\begin{equation} \label{2.4}
\frac{1}{2} \: \delta d (\omega, \omega) \: = \:
(\omega, \triangle \omega) \: - \: (\nabla \omega, \nabla \omega) 
\: - \: ( \omega, \Ric \: \omega).
\end{equation}
On the other hand,
\begin{equation} \label{2.5}
\frac{1}{2} \: i_{(d \ln \phi)_\sharp} d  (\omega, \omega) \:  = \:
\frac{1}{2} \: {\cal L}_{(d \ln \phi)_\sharp}  (\omega, \omega)
\end{equation}
We have
\begin{equation} \label{2.6}
{\cal L}_{(d \ln \phi)_\sharp} g \: = \: 2 \: \Hess(\ln \phi)
\end{equation}
Then
\begin{equation} \label{2.7}
\frac{1}{2} \: i_{(d \ln \phi)_\sharp} d  (\omega, \omega) \:  = \:
(\omega, {\cal L}_{(d \ln \phi)_\sharp} \omega) \: - \: 
(\omega, \Hess(\ln \phi) \: \omega).
\end{equation}
(The minus sign in (\ref{2.7}) comes from the fact that the pairing is
on $1$-forms instead of vector fields.)
Equations (\ref{2.2}), (\ref{2.3}), (\ref{2.4}) and (\ref{2.7}) give
\begin{equation} \label{2.8} 
\frac{1}{2} \: \widetilde{\delta} d (\omega, \omega) \: = \:
(\omega, \widetilde{\triangle} \omega) \: - \: (\nabla \omega, \nabla \omega) 
\: - \: (\omega, \BE \: \omega).
\end{equation}
Multiplying (\ref{2.8}) by $\phi$ and integrating over $M$, we obtain
\begin{equation} \label{2.9} 
0 \: = \: \langle \omega,  \widetilde{\triangle} \omega \rangle 
\: - \: 
\langle \nabla \omega, \nabla \omega \rangle 
\: - \: 
\langle \omega, \BE \:  \omega \rangle,
\end{equation}
or
\begin{equation} \label{2.10} 
\langle d \omega, \: d \omega \rangle \: + \: 
\langle \widetilde{\delta} \omega, \: \widetilde{\delta} \omega \rangle
 \: - \: \langle \nabla \omega, \: \nabla \omega \rangle
\: = \: \langle \omega, \:
\BE \: \omega \rangle.
\end{equation}

We can apply usual elliptic theory to the de Rham complex, with
the inner product  $\langle \cdot , \cdot \rangle$, to obtain an
isomorphism
\begin{equation} \label{2.11} 
\HH^*(M; \R) \: \cong \: \{ \omega \in \Omega^*(M) \: : \: 
d \omega \: = \: \widetilde{\delta} \omega \: = \: 0 \}. 
\end{equation}
If 
$\BE \: \ge \: 0$ and a $1$-form $\omega$ satisfies
$d\omega \: = \: \widetilde{\delta} \omega \: = \: 0$
then (\ref{2.10}) implies that
$\nabla \omega \: = \: 0$. Hence 
${\delta} \omega \: = \: 0$. Along with  
$\widetilde{\delta} \omega \: = \: 0$, (\ref{2.2}) now implies that
$\omega( \grad (\phi)) \: = \: 0$. Conversely, if
$\nabla \omega \: = \: \omega( \grad (\phi)) \: = \: 0$ then
$d \omega \: = \: \widetilde{\delta} \omega \: = \: 0$.
This proves the isomorphism in
part 3 of the theorem. 

If $\bb_1(M) \: = \: n$ then
there are $n$ linearly-independent
parallel $1$-forms on $M$ that annihilate
$\grad (\phi)$. The usual argument shows that $M$ is a flat torus.
As the parallel $1$-forms on $M$ annihilate
$\grad (\phi)$, $\phi$ must be constant. This proves part 3 of the
theorem.

A pointwise algebraic computation shows that
\begin{equation} \label{2.12} 
( d \omega, \: d \omega) \: + \: 
( {\cal L}_{\omega_\sharp} g, \: 
{\cal L}_{\omega_\sharp} g )
 \: = \: 2 \: ( \nabla \omega, \: \nabla \omega ).
\end{equation}

Then (\ref{2.10}) becomes
\begin{equation} \label{2.13} 
\langle \nabla \omega, \: \nabla \omega \rangle \: + \: 
\langle \widetilde{\delta} \omega, \: \widetilde{\delta} \omega \rangle
 \: - \: \langle \omega, \:
\BE \: \omega \rangle
\: = \: \langle {\cal L}_{\omega_\sharp} g, \: 
{\cal L}_{\omega_\sharp} g \rangle.
\end{equation}

If $\BE \: < \: 0$ and ${\cal L}_V g \: = \: 0$ then taking
$\omega \: = \: V^\sharp$,
(\ref{2.13}) implies that $V \: = \: 0$.  Hence the isometry group
of $(M, g)$ is discrete and, being compact, must be finite.
This proves part 4 of the theorem.

If $\BE \: \le \: 0$ and ${\cal L}_V g \: = \: 0$ then
(\ref{2.13}) implies that $\nabla V^\sharp \: = \:
\widetilde{\delta} V^\sharp \: = \: 0$. 
As before, we obtain that $V \phi \: = \: 0$. This proves part
5 of the theorem.\\ \\
{\bf Remarks : } \\
1. If we put $\omega \: = \: df$ in (\ref{2.8}) then we
recover the definition of $\BE$ from \cite{Bakry-Emery (1985)}. \\ \\
2. Jianguo Cao pointed out to me that a formula related to (\ref{2.10}) has 
been used to study the $\overline{\partial}$-operator on complete
K\"ahler manifolds \cite[Th\'eor\`eme 5.1]{Demailly (1982)}. \\ \\
3. The operator $\widetilde{\Delta}$ is related to the Witten Laplacian of
\cite{Witten (1982)}, but the two operators are distinct. To see
the relation, note that $\widetilde{\delta} \: = \: 
\phi^{- \: 1} \: \delta \:
{\phi}$. Put $D \: = \: \phi^{\frac12} \: d \:
\phi^{- \: \frac12}$ and $D^* \: = \: \phi^{- \: \frac12} \: \delta \:
\phi^{\frac12}$. Then the Witten Laplacian 
$D D^* \: + \: D^* D$ is related to
$\widetilde{\Delta}$ by 
\begin{equation}
D D^* \: + \: D^* D \: = \: \phi^{\frac12} \: \widetilde{\Delta} \:
\phi^{- \: \frac12}.
\end{equation}
The Bochner-type identity (\ref{2.10}), when translated to a statement about
$D D^* \: + \: D^* D$, becomes
\begin{equation}
D D^* \: + \: D^* D \: = \: \left( \phi^{\frac12} \: \nabla \:
\phi^{- \: \frac12} \right)^* \: \left( \phi^{\frac12} \: \nabla \:
\phi^{- \: \frac12} \right) \: + \: \BE,
\end{equation}
where the adjoints are with respect to the
unweighted $L^2$-inner product.
In contrast, in Morse-Witten theory one collects the terms differently,
by writing
$D D^* \: + \: D^* D \: = \: \nabla^* \nabla \: + \: \ldots$. \\ \\
4. The equality (\ref{2.10}) gives a way of defining the notion of
$\BE \: \ge \: r \: g$ for a class of nonsmooth measured manifolds
$(M, g, \phi)$. Namely, suppose that $M$ is a manifold whose
transition maps are $C^{1,1}$-regular. Let $g$ be a Riemannian
metric on $M$ whose components, in local charts, are in $C^0 \cap
H^1$, where $H^1$ denotes the Sobolev space. Let $\phi \in
C^0(M) \cap H^1_{loc}(M)$ be a positive function.  (There are a smooth
manifold $M^\prime$ and a $C^{1,1}$-diffeomorphism $M^\prime \rightarrow
M$. Hence after pulling back, if one wants then one can assume that $g$ 
and $\phi$ are defined on a smooth manifold.) 

\begin{definition} \label{definition1}
We say that $\Ric(M, g) \: \ge \: r \: g$ if for all compactly-supported
Lipschitz-regular $1$-forms $\omega$ on $M$,
\begin{equation} \label{2.15}
\int_M (d \omega, \: d \omega) \: d\vol_M \: + \: 
\int_M ({\delta} \omega, \: {\delta} \omega) \: d\vol_M
 \: - \: \int_M (\nabla \omega, \: \nabla \omega) \: d\vol_M
\: \ge \: r \: \int_M ( \omega, \: \omega) \: d\vol_M.
\end{equation}
We say that
$\BE (M, g, \phi) \: \ge \: r \: g$ if for all compactly-supported
Lipschitz-regular
$1$-forms $\omega$ on $M$,
\begin{equation} \label{2.16}
\langle d \omega, \: d \omega \rangle \: + \: 
\langle \widetilde{\delta} \omega, \: \widetilde{\delta} \omega \rangle
 \: - \: \langle \nabla \omega, \: \nabla \omega \rangle
\: \ge \: r \: \langle \omega, \: \omega \rangle.
\end{equation}
We say that
$\BEq (M, g, \phi) \: \ge \: r \: g$ if for all compactly-supported
Lipschitz-regular
$1$-forms $\omega$ on $M$,
\begin{equation} \label{2.17}
\langle d \omega, \: d \omega \rangle \: + \: 
\langle \widetilde{\delta} \omega, \: \widetilde{\delta} \omega \rangle
 \: - \: \langle \nabla \omega, \: \nabla \omega \rangle
\: - \: \frac{1}{q} \: \int_M \left( \omega(\nabla \ln \phi) \right)^2 \:
\phi \: d\vol_M
\: \ge \: r \: \langle \omega, \: \omega \rangle.
\end{equation}

\end{definition}

An immediate consequence of the definition is the following lemma.

\begin{lemma} \label{lemma1}
Let $M$ be a smooth  closed manifold. \\
1. If $\{g_i\}_{i=1}^\infty$
is a sequence of smooth Riemannian metrics on $M$
with $\Ric(M, g_i) \: \ge \: r \: g_i$, and 
$g_i \stackrel{C^0 \cap H^1}{\longrightarrow} g$ for some
$C^0 \cap H^1$-regular metric $g$, then
$\Ric(M, g) \: \ge \: r \: g$.\\
2.
If $\{(g_i, \phi_i)\}_{i=1}^\infty$
is a sequence of smooth Riemannian metrics and smooth
positive functions on $M$
with $\BE(M, g_i, \phi_i) \: \ge \: r \: g_i$, and 
$(g_i, \phi_i) \: 
\stackrel{C^0 \cap H^1}{\longrightarrow} \: (g, \phi)$ for some
$C^0 \cap H^1$-regular pair $(g, \phi)$, then
$\BE(M, g, \phi) \: \ge \: r \: g$.\\
3. If $\{(g_i, \phi_i)\}_{i=1}^\infty$
is a sequence of smooth Riemannian metrics and smooth
positive functions on $M$
with $\BEq(M, g_i, \phi_i) \: \ge \: r \: g_i$, and 
$(g_i, \phi_i) \: 
\stackrel{C^0 \cap H^1}{\longrightarrow} \: (g, \phi)$ for some
$C^0 \cap H^1$-regular pair $(g, \phi)$, then
$\BEq(M, g, \phi) \: \ge \: r \: g$.  
\end{lemma}

For example, let $\{ (M_i, g_i) \}_{i=1}^\infty$ be a sequence
of $n$-dimensional closed
Riemannian manifolds with Ricci curvatures bounded below by $r \in \R$,
injectivity radii bounded below by $i_0 \in \R^+$ and diameters bounded
above by $D \in \R^+$. Then $\{ (M_i, g_i) \}_{i=1}^\infty$ has a limit
point $X$ in the Gromov-Hausdorff topology. From
\cite{Anderson-Cheeger (1992)}, $X$ is an $n$-dimensional closed
manifold with a Riemannian 
metric $g$ that is $W^{1,p}$-regular for all 
$p \in [1, \infty)$. From the Sobolev embedding theorem, $g$ is also
$C^{0, \alpha}$-regular for all $\alpha \in (0,1)$.
After applying diffeomorphisms one has
$W^{1,p}$-convergence of a subsequence of
$\{ (M_i, g_i) \}_{i=1}^\infty$ to $(X, g)$, and so
$\Ric(X, g) \: \ge \: r \: g$ in the sense of Definition
\ref{definition1}. 

For another example, suppose that $M$ is a compact K\"ahler manifold
with local complex coordinates $\{ z^\alpha \}$ and metric
$g_{\alpha \overline{\beta}}$. Its Ricci form, in local
coordinates, is the $(1,1)$-form
$- \: \frac12 \: \partial \overline{\partial} \ln \det (g)$.
Now suppose that the $g_{\alpha \overline{\beta}}$ are only 
$C^0 \cap H^1$-regular.  The K\"ahler condition still makes sense
distributionally, and the Ricci form makes sense as a closed
$(1, 1)$-current.  Then $\Ric(M, g) \: \ge \: 0$ in the sense of
Definition \ref{definition1} if and only if 
$- \: \frac12 \: \partial \overline{\partial} \ln \det (g)$
is a positive current.
(This last condition makes sense for a much larger class of $g$.)

\section{Proof of Theorem \ref{theorem2}} \label{Section 3}

We (mostly) use the notation of \cite[Chapter 9]{Besse (1987)}.
If $X$ is a vector field on $B$, let $\overline{X}$ be its 
horizontal lift to $M$.
Let $N$ be the mean curvature vector field to the
fibers $F$. Let $A$ be the curvature of the horizontal distribution.
Let $T$ be the second fundamental
form tensor of the fibers $F$. Let $\nabla^M$ be the covariant derivative
operator on $M$ and let $\nabla^B$ be the covariant derivative operator on 
$B$.  From \cite[(9.36c)]{Besse (1987)}, there is an identity of
functions on $M$ :
\begin{equation} \label{3.1} 
\Ric^M( \overline{X}, \overline{X}) \: = \:
\Ric^B(X, X) \: - \: 2 \: \left( A_{\overline{X}}, 
A_{\overline{X}} \right) \: - \: \left( T \overline{X}, T
\overline{X} \right)
\: + \: \left( \overline{X}, \nabla^M_{\overline{X}} N \right). 
\end{equation}
  
Given $b \in B$,
let $\{\theta_t\}_{t \in (- \epsilon, 
\epsilon)}$ be the flow of ${X}$ as defined in a neighborhood of
$b$ and for $t$ in some interval $(- \epsilon, 
\epsilon)$. Let $\{\overline{\theta}_t\}_{t \in (- \epsilon, 
\epsilon)}$ be the flow of $\overline{X}$ that
covers $\theta_t$. It sends fibers to fibers diffeomorphically. Hence it 
makes sense to define
${\cal L}_{\overline{X}} \:  d\vol_F$ by 
\begin{equation} \label{3.2}
({\cal L}_{\overline{X}} \:  d\vol_F) \Big|_{F_b} \: = \:
\frac{d}{dt} \Big|_{t = 0} ({\overline{\theta}_t}^* d\vol_F) \Big|_{F_b}.
\end{equation}
With our conventions,
\begin{equation} \label{3.3} 
{\cal L}_{\overline{X}} \:  d\vol_F \: = \: - \: \left( 
\overline{X}, N \right) \:
d\vol_F.
\end{equation}
We have
\begin{equation} \label{3.4} 
\phi^B \: = \: \int_F \phi^M \: d\vol_F. 
\end{equation}
Then
\begin{align} \label{3.5}
X \phi^B \: & = \: {\cal L}_X \phi^B \: = \:
{\cal L}_X \int_F \phi^M \: d\vol_F  = \: 
\int_F {\cal L}_{\overline{X}} \: (\phi^M d\vol_F) \\
& = \: \int_F \left( \overline{X} \phi^M \: - \:
\left( \overline{X}, N \right) \phi^M \right) \: d\vol_F \notag 
\end{align}
and
\begin{align} \label{3.6} 
X X \phi^B \: & = \: \int_F
\left[\overline{X} \left( \overline{X} \phi^M \: - \:
\left( \overline{X}, N \right) \phi^M \right) \: - \:
\left( \overline{X}, N \right) \: 
\left( \overline{X} \phi^M \: - \:
\left( \overline{X}, N \right) \phi^M \right) \right]
 \: d\vol_F \\
& = \: \int_F
\left[ \overline{X}  \overline{X} \phi^M \: - \:
\overline{X} \left( \overline{X}, N \right) \: \phi^M  \: - \:
2 \: \left( \overline{X}, N \right) \: 
\overline{X} \phi^M \: + \:
\left( \overline{X}, N \right)^2 \phi^M  \right]
 \: d\vol_F \notag \\
& = \: \int_F
\left[ \frac{\overline{X}  \overline{X} \phi^M}{\phi^M} \:
- \: \left( \nabla^M_{\overline{X}} \overline{X}, N \right) \: 
- \: \left( \overline{X}, \nabla^M_{\overline{X}} N \right)
\: - \:
\left( \frac{\overline{X} \phi^M}{\phi^M} \right)^2 \right. \notag \\
& \: \: \: \: \: \: \: \:  \: \: \: \: \: \: \: \: \left. + \:
\left( \frac{\overline{X} \phi^M}{\phi^M} \: - \:
\left( \overline{X}, N \right) \right)^2  \right] \phi^M
 \: d\vol_F. \notag
\end{align}
Using the fact that $\nabla^M_{\overline{X}} \overline{X} \: = \:
\overline{\nabla^B_X X}$ \cite[(9.25d)]{Besse (1987)}, 
it follows that
\begin{align} \label{3.7} 
\Hess(\phi_B)(X, X) \: & = \:
X X \phi^B \: - \: (\nabla^B_X X) \phi^B \\
& = \: \int_F
\left[ \frac{\Hess(\phi^M)(\overline{X},\overline{X})}{\phi^M} \:
- \: \left( \overline{X}, \nabla^M_{\overline{X}} N \right)
\: - \:
\left( \frac{\overline{X} \phi^M}{\phi^M} \right)^2 \right. \notag \\
& \: \: \: \: \: \: \: \:  \: \: \: \: \: \: \: \: \left. + \:
\left( \frac{\overline{X} \phi^M}{\phi^M} \: - \:
\left( \overline{X}, N \right) \right)^2  \right] \phi^M
 \: d\vol_F \notag \\
& = \: \int_F
\left[ \Hess(\ln \phi^M)(\overline{X},\overline{X}) \:
- \: \left( \overline{X}, \nabla^M_{\overline{X}} N \right) \right. \notag \\
& \: \: \: \: \: \: \: \:  \: \: \: \: \: \: \: \: \left. + \:
\left( \frac{\overline{X} \phi^M}{\phi^M} \: - \:
\left( \overline{X}, N \right) \right)^2  \right] \phi^M
 \: d\vol_F. \notag
\end{align}
Substituting $\left( \overline{X}, \nabla^M_{\overline{X}} N \right)$
from (\ref{3.1}) gives
\begin{align} \label{3.8} 
\Ric^B(X, X) \: \phi^B \: - \: \Hess(\phi^B)(X, X) \:
& = \: \int_F
\left[ \BE^M( \overline{X}, \overline{X} ) \: + \: 
2 \: \left( A_{\overline{X}}, 
A_{\overline{X}} \right) \: + \: \left( T \overline{X}, T
\overline{X} \right) 
\right. \\
& \: \: \: \: \: \: \: \:  \: \: \: \: \: \: \: \: \left. 
- \: \left( \frac{\overline{X} \phi^M}{\phi^M} \: - \:
\left( \overline{X}, N \right) \right)^2  \right] \phi^M
 \: d\vol_F \notag
\end{align}
Using (\ref{3.5}),
\begin{align} \label{3.9} 
\BE^B(X, X) \: \phi^B \: 
& = \: \left[ \Ric^B(X, X)  \: - \:
\frac{\Hess(\phi^B)(X, X)}{\phi^B} \: + \:
\frac{(X \phi^B)^2}{(\phi^B)^2} \right] \: \phi^B \\ 
& = \: \int_F
\left[ \BE^M( \overline{X}, \overline{X} ) \: + \: 
2 \: \left( A_{\overline{X}}, 
A_{\overline{X}} \right) \: + \: \left( T \overline{X}, T
\overline{X} \right) 
\right. \notag \\
& \: \: \: \: \: \: \: \:  \: \: \: \: \: \: \: \: \left. 
- \: \left( \frac{\overline{X} \phi^M}{\phi^M} \: - \:
\left( \overline{X}, N \right) \right)^2  \right] \phi^M
 \: d\vol_F \notag \\
& + \: \left( \int_F \left( \frac{\overline{X} \phi^M}{\phi^M} \: - \:
\left( \overline{X}, N \right) \right) \: \phi^M
 \: d\vol_F \right)^2 \: (\phi^B)^{-1}. \notag
\end{align}

We have
\begin{equation} \label{3.10}
{\cal L}_{\overline{X}} (\phi^M d\vol_F) \: = \:
\left( \frac{\overline{X} \phi^M}{\phi^M} \: - \: (\overline{X}, N)
\right) \: \phi^M \: d\vol_F.
\end{equation}
By assumption, $\frac{\overline{X} \phi^M}{\phi^M} \: - \: 
(\overline{X}, N)$ is constant on a fiber $F$. Then
\begin{align} \label{3.11} 
\BE^B(X, X) \: \phi^B \: & = \: \int_F
\left[ \BE^M( \overline{X}, \overline{X} ) \: + \: 
2 \: \left( A_{\overline{X}}, 
A_{\overline{X}} \right) \: + \: \left( T \overline{X}, T
\overline{X} \right) 
\right] \: \phi^M
 \: d\vol_F \\
& \ge \: \int_F
\BE^M( \overline{X}, \overline{X} )
\: \phi^M
 \: d\vol_F. \notag
\end{align}
If $\BE^M( \overline{X}, \overline{X} ) \: \ge \: r \:
g^M( \overline{X}, \overline{X} )$ then (\ref{3.11}) implies that
$\BE^B( {X}, {X} ) \: \ge \: r \:
g^B( {X}, {X} )$. This proves Theorem \ref{theorem2}.1.

Now suppose that $\phi^M \: = \: 1$. Equations
(\ref{1.2}) and (\ref{3.9}) imply that  
\begin{align} \label{3.12} 
\BEq^B(X, X) \: \phi^B \:
 = \: & \int_F
\left[ \Ric^M( \overline{X}, \overline{X} ) \: + \: 
2 \: \left( A_{\overline{X}}, 
A_{\overline{X}} \right) \: + \: \left( T \overline{X}, T
\overline{X} \right) \:
- \: \frac{1}{q} \: \left( \overline{X}, N \right)^2  \right] 
 \: d\vol_F \\
& + \: \left( 1 \: - \: \frac{1}{q} \right) \: 
\left( - \: 
\int_F \left( \overline{X}, N \right)^2 \: d\vol_F \: + \: 
\left( \int_F \left( \overline{X}, N \right) 
 \: d\vol_F \right)^2  \: (\phi^B)^{-1}
\right). \notag
\end{align}
As $\left( \overline{X}, N \right) \: = \: - \: 
\Tr \left( T \overline{X} \right)$, we know that
$\left( T \overline{X}, T
\overline{X} \right) \:
- \: \frac{1}{q} \: \left( \overline{X}, N \right)^2 \: \ge \: 0$.
By assumption, $\left( \overline{X}, N \right)$ is constant on a fiber
$F$. 
Then
\begin{align} \label{3.13} 
\BEq^B(X, X) \: \phi^B \:
= \: & \int_F
\left[ \Ric^M( \overline{X}, \overline{X} ) \: + \: 
2 \: \left( A_{\overline{X}}, 
A_{\overline{X}} \right) \: + \: \left( T \overline{X}, T
\overline{X} \right) \:
- \: \frac{1}{q} \: \left( \overline{X}, N \right)^2  \right] 
 \: d\vol_F \\
\ge \: & \int_F
\Ric^M( \overline{X}, \overline{X} ) \: d\vol_F. \notag
\end{align}
If $\BE^M( \overline{X}, \overline{X} ) \: \ge \: r \:
g^M( \overline{X}, \overline{X} )$ then
\begin{equation} \label{3.14}
\BEq^B(X, X) \: \phi^B \:
\ge \: r \: \int_F
g^M( \overline{X}, \overline{X} ) \: d\vol_F \: = \:
r \: g^B(X, X) \: \phi^B. 
\end{equation}
This proves Theorem \ref{theorem2}.2. \\ \\
{\bf Example : } Let $p \: : \: M \rightarrow B$ be a Riemannian
submersion, with $M$ compact,
whose fiber transport preserves the fiberwise metric up to
multiplicative constants.  Equivalently, the Riemannian metric $g$
on $M$ comes from
starting with a submersion
metric $g^\prime$ with totally geodesic fibers, along with a
positive function $f \in C^\infty(B)$, and then
multiplying the fiberwise metric of $g^\prime$ on $F_b$ by $f^2(b)$. 
One can think of $g$ as a generalized warped product metric.

Suppose that the fibers $F$ have nonnegative Ricci curvature. For 
$\epsilon \: > \: 0$, let $g_\epsilon$ be the Riemannian metric on $M$ 
which comes
from multiplying the fiberwise Riemannian metrics by $\epsilon^2$. 
Then as $\epsilon \rightarrow 0$, the metrics $g_\epsilon$ have Ricci
curvatures that are uniformly bounded below. Explicitly, let 
$\overline{X}$ be the
horizontal lift of a vector field $X$ on $B$ and let 
$\overline{U}$ be a vertical vector field.
Then as $\epsilon \rightarrow 0$, with the notation of
\cite[Chapter 9]{Besse (1987)},
\begin{align} \label{3.15}
\Ric^M_\epsilon(\overline{X}, \overline{X}) & \sim
p^* \: \Ric^B({X}, {X}) \: - \:
(T\overline{X}, T\overline{X}) \: + \: 
\left( \overline{X}, \nabla^M_{\overline{X}} N \right), \\
\Ric^M_\epsilon(\overline{X}, \overline{U}) & \sim 0 \notag \\
\Ric^M_\epsilon(\overline{U}, \overline{U}) & \sim
\Ric^F(\overline{U}, \overline{U}) \: + \:
\epsilon^2 \left( (\widetilde{\delta} T)(\overline{U}, \overline{U})
\: - \: (N, T_{\overline{U}} \overline{U}) \right). \notag
\end{align}
(The terms on the right-hand-side of (\ref{3.15}) are evaluated with
respect to the
metric $g_1$.) This is an example of a collapse with Ricci curvature
bounded below, to which Theorem \ref{theorem2}.2 applies. 

For another example, let $M$ be a compact Riemannian manifold on which
a Lie group $G$ acts isometrically and effectively. Suppose that the
$G$-action on $M$ has a single orbit type and put $B \: = \: G \backslash M$.
Then there is a natural Riemannian submersion $p \: : \: M \rightarrow B$.
As the orbits of the $G$-action on $M$ are all $G$-diffeomorphic to a 
homogeneous space $G/H$, and $G/H$ has a unique $G$-invariant volume form
up to constants, it follows that the fiber transport of the Riemannian
submersion preserves measures up to constants.  Hence Theorem
\ref{theorem2}.2 applies.

\section{Proof of Theorem \ref{theorem3}} \label{Section 4}

We refer to \cite{Fukaya (1987)} for the definition of the measured
Gromov-Hausdorff topology.

To prove Theorem \ref{theorem3}.1, we just apply
the warped product construction of the proof of
Theorem \ref{theorem1}.1 to $S^q \times B$.

Let $\{M_i, g_i \}_{i=1}^\infty$ be a sequence as in
the statement of Theorem \ref{theorem3}.2. We may assume that
$\lim_{i \rightarrow \infty} (M_i, g_i, d\vol_i) \: = \:
(X, \mu)$ in the measured Gromov-Hausdorff topology.
If $q \: = \: 0$ then $X$ is a smooth manifold with a $C^{1,\alpha}$-regular
metric $g^X$ and after taking a subsequence and applying diffeomorphisms, we
may assume that $(M_i, g_i)$ converges to $(X, g^X)$ in the
$C^{1,\alpha}$-topology
(see, for example, \cite{Kasue (1989)}). In this case, the theorem
follows from Lemma \ref{lemma1}.1.

Suppose that $q \: > \: 0$.
By saying that $X$ is a manifold, we mean that in the construction of
$X$ as a quotient space $\widehat{X}/O(N)$ \cite{Fukaya (1988)}, the
action of $O(N)$ on the manifold $\widehat{X}$ has a single orbit type.
Then $X$ has the structure of a smooth manifold with a
$C^{1, \alpha}$-regular pair $(g^X, \phi^X)$.

For any $\epsilon \: > \: 0$, we
can apply smoothing results of Abresch and others 
\cite[Theorem 1.12]{Cheeger-Fukaya-Gromov (1992)}
to obtain new metrics $g_i(\epsilon)$ with
\begin{align} \label{4.1}
e^{- \: \epsilon} \: g_i \: \le \: g_i(\epsilon) \: & \le \: 
e^{\epsilon} \: g_i, \\
|\nabla_{g_i} \: - \: \nabla_{g_i(\epsilon)} | \: & \le \: 
\epsilon, \notag \\
|\nabla^k_{g_i(\epsilon)} 
\: \Riem(M_i, g_i(\epsilon))|
\: & \le \: C_{k}(N, \epsilon, \Lambda), \notag
\end{align}
where the constants are uniform.
We can also assume that
$\Ric(M_i, g_i(\epsilon)) \: \ge \: (r \: - \: \epsilon) \: 
g_i(\epsilon)$ \cite[Remark 2, p. 51]{Dai-Wei-Ye (1996)}. (See
\cite[Theorem 2.1]{Rong (1996)} for a similar statement about
sectional curvature.) 
For small $\epsilon$,
let $B(\epsilon)$ be a Gromov-Hausdorff limit of a subsequence 
of $\{ (M_i, g_i(\epsilon)) \}_{i=1}^\infty$.
We relabel the subsequence as $\{ (M_i, g_i(\epsilon)) \}_{i=1}^\infty$.
From \cite[Proposition 4.9]{Cheeger-Fukaya-Gromov (1992)},
for large $i$, there is a
small $C^2$-perturbation $g_i^\prime(\epsilon)$ of $g_i(\epsilon)$
which is invariant with respect to a $Nil$-structure.
In particular,
we may assume that $\Ric(M_i, g_i^\prime(\epsilon)) \: \ge \: 
(r \: - \: 2 \: \epsilon) \: g_i^\prime(\epsilon)$. Now 
$(M_i, g_i^\prime(\epsilon))$ is the
total space of a Riemannian submersion $M_i \rightarrow B(\epsilon)$ with 
infranil fibers and affine holonomy.
Let $\left( g_i^{B(\epsilon)}, \phi_i^{B(\epsilon)} \right)$ denote the
induced metric and measure on $B(\epsilon)$. As the fiber transport of the
Riemannian submersion preserves the affine-parallel volume forms of the fibers,
up to constants, Theorem \ref{theorem2}.2 implies that
$\BEq \left(
B(\epsilon), g_i^{B(\epsilon)}, \phi_i^{B(\epsilon)}
\right) \: \ge \: 
(r \: - \: 2 \: \epsilon) \:
g_i^{B(\epsilon)}$. 
Varying $i$ and $\epsilon$, we can extract a subsequence of
$\left\{ \left(
B(\epsilon), g_i^{B(\epsilon)}, \phi_i^{B(\epsilon)}
\right) \right\}$ with 
$i \rightarrow
\infty$ and $\epsilon \rightarrow 0$ that converges in
the $C^{1, \alpha}$-topology to $(X, g^X, \phi^X)$. The theorem now
follows from Lemma \ref{lemma1}.3.

\section{Proof of Theorem \ref{theorem4}} \label{Section 5}

Let $s$ be a segment from $t_0 \in T_0$ to $t \in T$, with length
$l(s) \: > \: u_3$ and arc-length
parameter $u$.  By definition, $s$ is length-minimizing.
We can decompose the measure $\phi \: d\vol_M$ on 
$A(u_1, u_4)$ as $\phi \: \area_s(u) \: 
du \: \mu(s)$, where $\mu$ is a measure on the space
${\cal S}$ of distinct segments $s$ that make up $A(u_1, u_4)$,
$du$ is the length measure along a segment $s$ and
$\area_s(u)$ is the relative size of the transverse
Riemannian area density
along $s$, as measured with respect to the fan of segments.
Let $h$ denote the trace of the second fundamental form $\Pi$ of a level
set of constant distance from $T_0$.
(With our conventions, the boundary of the unit ball in $\R^n$ has
positive mean curvature.)
Differentiating along $s$, with respect to $u$, gives
\begin{equation} \label{5.1}
\partial_u  \ln(\phi(u) \: \area_s(u)) \: \equiv \:
\frac{\partial_u (\phi(u) \: \area_s(u))}{\phi(u) \: 
\area(u))}
\: = \:
h(u) \: + \: \partial_u \ln \phi(u)
\end{equation}
and
\begin{equation} \label{5.2}
\partial^2_u  \ln(\phi(u) \: \area_s(u)) \: = \:
\partial_u h(u) \: + \: \partial^2_u \ln \phi(u).
\end{equation}
From the Riccati equation for $\Pi$,
\begin{equation} \label{5.3}
\partial_u h(u) \: = \: - \: \Tr (\Pi^2) \: - \: \Ric(\partial_u, \partial_u)
\: \le \: - \: \Ric(\partial_u, \partial_u).
\end{equation}
Then
\begin{equation} \label{5.4}
\partial^2_u  \ln(\phi(u) \: \area_s(u)) \: \le \:
- \: \BE(\partial_u, \partial_u) \: \le \: - \: r.
\end{equation}
Hence for any $c \in \R$,
\begin{equation} \label{5.5}
\partial^2_u  \left( 
\ln \left( \phi(u) \: \area_s(u)
\right) \: + \:
\frac{r}{2} \: u^2 \: - \: c \: u \right) \: \le \: 0.
\end{equation}

Fix $s$ and put
\begin{equation} \label{5.6}
a(u) \: = \: \phi(u) \: \area_s(u),
\end{equation} 
\begin{equation} \label{5.7}
\widehat{a}(u) \: = \: e^{- \: \frac{r}{2} \: u^2 \: + \: c \: u},
\end{equation} 
\begin{equation} \label{5.8}
v(u_1, u_2) \: = \:
\int_{u_1}^{u_2} a(u) \: du
\end{equation}
and 
\begin{equation} \label{5.9}
\widehat{v}(u_1, u_2) \: = \:
\int_{u_1}^{u_2} \widehat{a}(u) \: du.
\end{equation}
Then (\ref{5.5}) says that 
\begin{equation} \label{5.10}
\frac{d^2}{du^2} \: \ln \left( \frac{a}{\widehat{a}} \right) \: \le \: 0,
\end{equation}
i.e. that $\ln \left( \frac{a}{\widehat{a}} \right)$ is concave in $u$.

\begin{lemma} \label{lemma2}
If $\frac{v(u_2, u_3)}{\widehat{v}(u_2, u_3)} \: \le \:
\frac{v(u_1, u_2)}{\widehat{v}(u_1, u_2)}$ then 
$\frac{a(u_3)}{\widehat{a}(u_3)} \: \le \:
\frac{v(u_2, u_3)}{\widehat{v}(u_2, u_3)}$.
\end{lemma}
\begin{pf}
Suppose that 
\begin{equation} \label{5.11}
\frac{a(u_3)}{\widehat{a}(u_3)} \: > \:
\frac{v(u_2, u_3)}{\widehat{v}(u_2, u_3)} \: = \:
\frac{\int_{u_2}^{u_3} \frac{a(u)}{\widehat{a}(u)} \: 
\widehat{a}(u) \: du}{\int_{u_2}^{u_3}
\widehat{a}(u) \: du}. 
\end{equation}
If $\frac{a(u_2)}{\widehat{a}(u_2)} \: \ge \:
\frac{a(u_3)}{\widehat{a}(u_3)}$ then the concavity of 
$\ln \left( \frac{a}{\widehat{a}} \right)$ implies that
 \begin{equation} \label{5.12}
\frac{a(u_3)}{\widehat{a}(u_3)} \: \le \:
\frac{\int_{u_2}^{u_3} \frac{a(u)}{\widehat{a}(u)} \: 
\widehat{a}(u) \: du}{\int_{u_2}^{u_3}
\widehat{a}(u) \: du}, 
\end{equation}
which is a contradiction.  Thus 
\begin{equation} \label{5.13}
\frac{a(u_2)}{\widehat{a}(u_2)} \: < \:
\frac{a(u_3)}{\widehat{a}(u_3)}.
\end{equation} 
With the concavity
of $\ln \left( \frac{a}{\widehat{a}} \right)$, (\ref{5.13}) implies that
$\frac{a(u)}{\widehat{a}(u)}$ is decreasing in $u$
for $u \: < \: u_2$ and so
\begin{equation} \label{5.14}
\frac{\int_{u_1}^{u_2} \frac{a(u)}{\widehat{a}(u)} \: 
\widehat{a}(u) \: du}{\int_{u_1}^{u_2}
\widehat{a}(u) \: du} \: < \:
\frac{a(u_2)}{\widehat{a}(u_2)}.
\end{equation}
The concavity
of $\ln \left( \frac{a}{\widehat{a}} \right)$
and (\ref{5.13}) also imply that
\begin{equation} \label{5.15}
\frac{a(u_2)}{\widehat{a}(u_2)} \: < \:
\frac{\int_{u_2}^{u_3} \frac{a(u)}{\widehat{a}(u)} \: 
\widehat{a}(u) \: du}{\int_{u_2}^{u_3}
\widehat{a}(u) \: du}.
\end{equation}
Thus we have
\begin{equation} \label{5.16}
\frac{\int_{u_1}^{u_2} \frac{a(u)}{\widehat{a}(u)} \: 
\widehat{a}(u) \: du}{\int_{u_1}^{u_2}
\widehat{a}(u) \: du} \: < \:
\frac{a(u_2)}{\widehat{a}(u_2)} \: < \:
\frac{\int_{u_2}^{u_3} \frac{a(u)}{\widehat{a}(u)} \: 
\widehat{a}(u) \: du}{\int_{u_2}^{u_3}
\widehat{a}(u) \: du},
\end{equation}
which contradicts the assumption.
\end{pf}

\begin{lemma} \label{lemma3}
If $\frac{v(u_2, u_3)}{\widehat{v}(u_2, u_3)} \: \le \:
\frac{v(u_1, u_2)}{\widehat{v}(u_1, u_2)}$ then for $u_4 \in (u_3, l(s))$,
$\frac{v(u_3, u_4)}{\widehat{v}(u_3, u_4)} \: \le \:
\frac{v(u_2, u_3)}{\widehat{v}(u_2, u_3)}$.
\end{lemma}
\begin{pf}
For $u \in (u_3, l(s))$, put
\begin{equation} \label{5.17}
F(u) \: = \: \ln \left( 
\frac{v(u_3, u)}{\widehat{v}(u_3, u)} \Big/
\frac{v(u_2, u_3)}{\widehat{v}(u_2, u_3)}
\right).
\end{equation}
Then 
\begin{equation} \label{5.18}
F^\prime(u) \: = \:
\frac{a(u)}{v(u_3, u)} \: - \: 
\frac{\widehat{a}(u)}{\widehat{v}(u_3, u)} \: = \:
\frac{\widehat{a}(u)}{v(u_3, u)} \:
\left[ \frac{a(u)}{\widehat{a}(u)} \: - \: 
\frac{v(u_3, u)}{\widehat{v}(u_3, u)} \right].
\end{equation}
Lemma \ref{lemma2} implies that if $F(u) \: \le \: 0$ then
$F^\prime(u) \: \le \: 0$. We can extend $F(u)$ smoothly to $u \: = \:
u_3$, with
\begin{equation} \label{5.19}
F(u_3) \: = \: \ln \left( 
\frac{a(u_3)}{\widehat{a}(u_3)} \Big/
\frac{v(u_2, u_3)}{\widehat{v}(u_2, u_3)}
\right).
\end{equation}
By Lemma \ref{lemma2}, $F(u_3) \: \le \: 0$. It follows that
$F(u) \: \le \: 0$ for all $u \: \in (u_3, l(s))$, which proves the lemma.
\end{pf}

We have
\begin{equation} \label{5.20}
\frac{\vol_\phi(A(u_2, u_3))}{\vol_\phi(A(u_1, u_2))} \: = \:
\frac{\int_{\cal S} \frac{v_s(u_2, u_3)}{v_s(u_1, u_2)} \: v_s(u_1, u_2) \: 
d\mu(s)}{\int_{\cal S} 
v_s(u_1, u_2) \: d\mu(s)}.
\end{equation}
Put 
\begin{equation} \label{5.21}
{\cal S}^\prime \: = \: \left\{s \in S \: : \: 
\frac{v_s(u_2, u_3)}{v_s(u_1, u_2)} \: < \: 
\frac{\widehat{v}(u_2, u_3)}{\widehat{v}(u_1, u_2)} \right\}
\end{equation}
and
\begin{equation} \label{5.22}
T^\prime \: = \: \bigcup_{s \in S^\prime} s.
\end{equation}
We claim that (\ref{1.7}) is satisfied.  If it is not satisfied, put
${\cal S}^{\prime \prime} \: = \: {\cal S} \: - \: {\cal S}^\prime$ and
$T^{\prime \prime} \: = \: T \: - \: T^\prime$. Then
\begin{equation} \label{5.23}
\frac{\vol_\phi(T^{\prime \prime} \cap A(u_1, u_2))}{\vol_\phi(A(u_1, u_2))}
 \: > \: 
\frac{\vol_\phi(A(u_2, u_3))}{\vol_\phi(A(u_1, u_2))} \:
\left( 
\frac{\widehat{v}(u_2, u_3)}{\widehat{v}(u_1, u_2)} \right)^{-1}.
\end{equation}
However, from the definition of $T^{\prime \prime}$,
\begin{align} \label{5.24}
\vol_\phi(A(u_2, u_3)) \: & \ge \:
\vol_\phi(T^{\prime \prime} \cap A(u_2, u_3)) \: = \:
\int_{{\cal S}^{\prime \prime}} 
\frac{v_s(u_2, u_3)}{v_s(u_1, u_2)} \: v_s(u_1, u_2) \: 
d\mu(s) \\
& \ge \: \int_{{\cal S}^{\prime \prime}} 
\frac{\widehat{v}(u_2, u_3)}{\widehat{v}(u_1, u_2)} \: v_s(u_1, u_2) \: 
d\mu(s) \: = \: \frac{\widehat{v}(u_2, u_3)}{\widehat{v}(u_1, u_2)}
\: \vol_\phi(T^{\prime \prime} \cap A(u_1, u_2)), \notag
\end{align}
which contradicts (\ref{5.23}).

If there is a cutpoint along $s$, with respect to its basepoint in $T_0$,
at $u_c \in (u_3, u_4)$ then we put
$v_s(u_3, u_4) \: = \: \int_{u_3}^{u_c} a_s(u) \: du$, and otherwise we
put $v_s(u_3, u_4) \: = \: \int_{u_3}^{u_4} a_s(u) \: du$.
Using Lemma \ref{lemma3},
\begin{equation}
\frac{\vol_\phi(T^\prime \cap A(u_3, u_4))}{\vol_\phi(
T^\prime \cap A(u_2, u_3))} \: = \:
\frac{\int_{{\cal S}^\prime} \frac{v_s(u_3, u_4)}{v_s(u_2, u_3)} \:
v_s(u_2, u_3) \:  
d\mu(s)}{\int_{{\cal S}^\prime} 
v_s(u_2, u_3) \: d\mu(s)} \: \le \:
\frac{\widehat{v}_s(u_3, u_4)}{\widehat{v}_s(u_2, u_3)}.
\end{equation}
This proves the first part of the theorem.

Suppose that there is a number $r \in \R$ so that for each tube $T$ and
$c \in \R$ satisfying (\ref{1.6}), there is a subtube $T^\prime$ satisfying
the properties of the theorem.  Given $m \in M$ and a unit vector
$v \in T_mM$, let $T_0$ be a hypersurface passing through
$m$ such that $T_m(T_0) \: = \: v^\perp$ and the second fundamental form
of $T_0$ at $m$ vanishes. Let $s$ be a 
minimizing segment with $s(0) \: = \: m$ and $s^\prime(0) \: = \: v$.
From (\ref{5.1}), 
\begin{equation} \label{5.25}
\frac{d}{du} \Big|_{u=0} (\ln (\phi(u) \: \area(u))
\: = \: v (\ln \phi).
\end{equation}
From (\ref{5.2}) and the Riccati equation,
\begin{equation} \label{5.26}
\frac{d^2}{du^2} \Big|_{u=0} (\ln (\phi(u) \: \area(u))
\: = \: - \: \BE(v, v).
\end{equation}

Put $c_0 \: = \: v (\ln \phi)$ and $r_0 \: = \: \BE(v, v)$.
Then for small $u$, 
\begin{equation}
\ln (\phi(u) \: \area(u)) \: \sim \:
\const \: + \: c_0 \: u \: - \: \frac{r_0}{2} \: u^2.
\end{equation}
For small $u_1 \: < \: u_2 \: < \: u_3 \: < \:  u_4$, we have
\begin{equation}
\frac{v(u_2, u_3)}{v(u_1, u_2)} \: \sim \:
\frac{\int_{u_2}^{u_3} e^{- \: \frac{r_0}{2} \: u^2 \: + \: c_0 \: u \:
} du}{\int_{u_1}^{u_2} e^{- \: \frac{r_0}{2} \: u^2 \: + \: c_0 \: u} du}
\end{equation}
and
\begin{equation}
\frac{v(u_3, u_4)}{v(u_2, u_3)} \: \sim \:
\frac{\int_{u_3}^{u_4} e^{- \: \frac{r_0}{2} \: u^2 \: + \: c_0 \: u}
du}{\int_{u_2}^{u_3} 
e^{- \: \frac{r_0}{2} \: u^2 \: + \: c_0 \: u} du}.
\end{equation}
Take $T$ to be a small tube around $s$ (with small base $T_0$), take
$u_3$ small relative to $u_4$ and take
$c \: = \: c_0 \: + \: \epsilon$ with $\epsilon \: > \: 0$ small so
that (\ref{1.6}) holds. If there is to be a subtube $T^\prime$ such that
(\ref{1.8}) holds, for all such choices, then
we must have $r_0 \: \ge \: r$. This proves the theorem.

\section{Remarks} \label{Section 6}
\noindent
1. If $M^n$ is compact and 
$\widetilde{\Ric}_q \: \ge \: r \: g$, with $q$ an integer greater
than one, then Theorem \ref{theorem3}.1 says that
$(M, g, \phi)$ is the limit
of a sequence of $(n \: + \: q)$-dimensional manifolds with
Ricci curvature bounded below by $r$. As in the proof of Theorem
\ref{theorem1}.2, we can then apply standard results
about manifolds with Ricci curvature bounded below, in order to obtain
conclusions about $(M, g, \phi)$.
For example, applying the Bishop-Gromov
inequality to the $(n \: + \: q)$-dimensional manifolds and
taking the limit, we
obtain a Bishop-Gromov-type inequality for the measures of the
distance balls in $M$.
Namely, let $\vol_\phi$ denote the weighted measure.
Then for $0 \: < \: u_1 \: < \: u_2$, 
$\frac{\vol_\phi(B_{u_2})}{\vol_\phi(B_{u_1})}$ is less than or
equal to the corresponding quantity in the $(n \: + \: q)$-dimensional
space form of Ricci curvature $r$.
If $r \: > \: 0$ then
applying Myers' theorem to the
$(n \: + \: q)$-dimensional manifolds and
taking the limit, we obtain that
$\diam(M) \: \le \: \pi \:  \sqrt{\frac{n \: + \: q \: - \: 1}{r} }$.
This gives alternative proofs of some results of 
Qian \cite[Corollary  2 and Theorem 5]{Qian (1997)} 
in the special case when $q$ is an integer
greater than one.  (The results of \cite{Qian (1997)} are valid for
all positive $q$.) One can also show that if $\BEq \: \ge \: r \: g$
with $q \in (0, \infty)$ then
$(M, g, \phi)$ satisfies the directional Bishop-Gromov inequality of
\cite[(A.2.2)]{Cheeger-Colding (1997)} with respect to a model space
of formal dimension $n \: + \: q$.
\\ \\
2. Similarly, if $q$ is an integer greater than one then
there are Sobolev inequalities for the 
$(n \: + \: q)$-dimensional collapsing manifolds 
\cite[Theorem 3, p. 397]{Berard (1988)}. Applying these inequalities
to functions that pullback from $M$,  we
obtain weighted Sobolev inequalities for
$M$. Namely, put $V \: = \: \int_M \phi \: d\vol_M$. Given
$\alpha, \beta \in [1, \infty)$ such that
$\alpha \: \le \frac{(n+q) \beta}{n+q-\beta}$,
let $\Sigma(n+q; \alpha, \beta)$ be the
Sobolev constant of the
standard $(n+q)$-sphere $S^{n+q}$, defined by
\begin{equation} \label{6.1}
\Sigma(n+q; \alpha, \beta) \: = \: \sup \left\{
\frac{\parallel f \parallel_\alpha}{\parallel df \parallel_\beta} \: : \:
f \in W^{1,\beta}(S^{n+q}), \: f \neq \: 0, \: \int_{S^{n+q}} f \: = \: 0
\right\}.
\end{equation}
Then if $\widetilde{\Ric}_q(M, g, \phi) 
\: \ge \: \frac{n+q-1}{R^2} \: g$, we have
\begin{align} \label{6.2}
\left( \int_M f^\alpha \: \phi \: d\vol_M \right)^{\frac{1}{\alpha}} \: \le \: 
& \Sigma(n+q; \alpha, \beta) \: R \: \left( 
\frac{V}{\vol(S^{n+q})} 
\right)^{{\frac{1}{\alpha}} \: - \: {\frac{1}{\beta}}} \:
\left( \int_M |\nabla f|^\beta \: \phi \: d\vol_M \right)^{\frac{1}{\beta}} 
\: + \\
& V^{{\frac{1}{\alpha}} \: - \: {\frac{1}{\beta}}} \: 
\left( \int_M f^\beta \: \phi \: d\vol_M \right)^{\frac{1}{\beta}} \notag
\end{align}
for $f \in W^{1, \beta}(M)$.
In the case $\beta \: = \: 2$, these
inequalities appeared in \cite{Bakry (1994)}. \\ \\
3. From the Bishop-Gromov-type inequalities, one can easily show that
for any $q, D \in \R^+$ and $r \in \R$, the space of Riemannian manifolds 
$(M, g)$ with a smooth positive probability measure $\phi \: d\vol_M$
satisfying $\widetilde{\Ric}_q(M, g, \phi) \: \ge \:
r \: g$ and $\diam(M, g) \: \le \: D$, taken modulo diffeomorphisms,
is precompact in the measured Gromov-Hausdorff topology. 

Since the relative volume in 
$\R^{n+q}$ of $B_{u_2}$ and $B_{u_1}$
is $\left( \frac{u_2}{u_1} \right)^{n+q}$, we cannot
expect any Bishop-Gromov-type comparison theorem for the masses of
balls in spaces with
$\BE$ bounded below, i.e. when $q \rightarrow \infty$ in
$\widetilde{\Ric}_q$. However, it is interesting that spaces with
$\BE \: \ge \: r \: g$ for $r \: > \: 0$ do admit isoperimetric
inequalities \cite{Bakry-Ledoux (1996)}. \\ \\
4. It is an interesting question whether there is a good synthetic
notion of a metric-measure space with Ricci curvature bounded below, 
in analogy to the notion of an Alexandrov space with curvature bounded
below.  See \cite[Appendix 2]{Cheeger-Colding (1997)} for discussion.
It is clear from Theorem \ref{theorem3}.1 that triples $(M, g, \phi)$
with $\BEq \: \ge \: r \: g$ are examples of metric-measure spaces with
generalized Ricci curvature bounded below by $r$, at least if $q$ is an
integer greater than one. 

There are various ways that one could
try to extend the notion of Ricci curvature
bounded below, from smooth metric-measure
spaces to more general metric-measure spaces. One could fix 
$q \in (0, \infty)$ and try to extend
the notion of having $\BEq \: \ge \: r \: g$.
Or one could consider all $q$ simultaneously, and say in particular that
a triple $(M, g, \phi)$
has generalized Ricci curvature bounded below by $r$ if
$\BEq \: \ge \: r \: g$ for some $q \in (0, \infty)$.
Or one could consider a triple $(M, g, \phi)$ to have 
generalized Ricci curvature bounded below by $r$ if 
$\BE \: \ge \: r \: g$. 

We note that there is a difference between
having $\BEq \: \ge \: r \: g$ for some $q \in (0, \infty)$
and having $\BE \: \ge \: r \: g$.
For example, if $r \: > \: 0$ and
$\BEq \: \ge \: r \: g$ for some
$q \in (0, \infty)$ then $M$ is compact \cite[Theorem 5]{Qian (1997)},
whereas if $\BE \: \ge \: r \: g$ then $M$ can be noncompact
(as in the case of $\R$ with $\phi(x) \: = \:
e^{- \: \frac{r}{2} \: x^2}$.) It is also easy to see that triples 
$(M, g, \phi)$ with
$\BE \: \ge \: 0$ generally do not satisfy the splitting principle.

If one does consider a triple $(M, g, \phi)$ with
$\BE \: \ge \: r \: g$ to be an admissible space with 
generalized Ricci curvature bounded below by $r$ then one has
a large class of examples.  For instance, from this
viewpoint it would be reasonable to say that
flat $\R^n$ with the measure $e^{-V} \: dx_1 \ldots dx_n$ has
nonnegative generalized Ricci curvature if $V$ is any convex function
on $\R^n$.

\end{document}